\documentclass[11pt]{amsart}
\usepackage{amssymb}

\allowdisplaybreaks

\begin{document}
\theoremstyle{plain}
\newtheorem{thm}{Theorem}[section]
\newtheorem{theorem}[thm]{Theorem}
\newtheorem{lemma}[thm]{Lemma}
\newtheorem{corollary}[thm]{Corollary}
\newtheorem{proposition}[thm]{Proposition}
%%%%%%%%%%%%%%%%%%%% Text roman %%%%%%%%%%%%%%%%%%%%%%%%%%%%%
\theoremstyle{definition}
\newtheorem{construction}[thm]{Construction}
\newtheorem{notations}[thm]{Notations}
\newtheorem{question}[thm]{Question}
\newtheorem{problem}[thm]{Problem}
\newtheorem{remark}[thm]{Remark}
\newtheorem{remarks}[thm]{Remarks}
\newtheorem{definition}[thm]{Definition}
\newtheorem{claim}[thm]{Claim}
\newtheorem{assumption}[thm]{Assumption}
\newtheorem{assumptions}[thm]{Assumptions}
\newtheorem{properties}[thm]{Properties}
\newtheorem{example}[thm]{Example}
\newtheorem{comments}[thm]{Comments}
\newtheorem{blank}[thm]{}
\newtheorem{observation}[thm]{Observation}
\newtheorem{defn-thm}[thm]{Definition-Theorem}
\newtheorem{conjecture}[thm]{Conjecture}

\newcommand{\sM}{{\mathcal M}}

%%%%%%%%%%%%%%%%%%%%%%%%%%%%%%%%%%%%%%%%%%%%%%%%%%%%%%%%%%%%%%

\title[Intersection Numbers and Automorphisms of Stable Curves]{Intersection Numbers and\\ Automorphisms of Stable Curves}


%    author one information
\author{Kefeng Liu}
\address{Center of Mathematical Sciences, Zhejiang University, Hangzhou, Zhejiang 310027, China;
                Department of Mathematics,University of California at Los Angeles,
                Los Angeles, CA 90095-1555, USA}
\curraddr{} \email{liu@math.ucla.edu, liu@cms.zju.edu.cn}
\thanks{}

%    author two information
\author{Hao Xu}
\address{Center of Mathematical Sciences, Zhejiang University, Hangzhou, Zhejiang 310027, China}
\curraddr{} \email{haoxu@cms.zju.edu.cn}
\thanks{}

%    \subjclass is required.
%\subjclass[2000]{Primary 14A20}

\date{}

\dedicatory{}

\maketitle

\section{Introduction}

Denote by $\overline{\sM}_{g,n}$ the moduli space of stable
$n$-pointed genus $g$ complex algebraic curves. We have the morphism
that forgets the last marked point
$$
\pi: \overline{\sM}_{g,n+1}\longrightarrow \overline{\sM}_{g,n}.
$$
Denote by $\sigma_1,\dots,\sigma_n$ the canonical sections of $\pi$,
and by $D_1,\dots,D_n$ the corresponding divisors in
$\overline{\sM}_{g,n+1}$. Let $\omega_{\pi}$ be the relative
dualizing sheaf, we have the following tautological classes on
moduli spaces of curves.
\begin{align*}
\psi_i&=c_1(\sigma_i^*(\omega_{\pi}))\\
\kappa_i&=\pi_*\left(c_1\left(\omega_{\pi}\left(\sum D_i\right)\right)^{i+1}\right)\\
\lambda_l&=c_l(\pi_*(\omega_{\pi})),\quad 1\leq l\leq g.
\end{align*}
The classes $\kappa_i$ were first introduced by Mumford \cite{Mu} on
$\overline{\sM}_g$, their generalization to $\overline{\sM}_{g,n}$
here is due to Arbarello-Cornalba \cite{Ar-Co}.

For background materials about the intersection theory of moduli spaces of curves, we refer to the book \cite{La-Zv}
and the survey paper \cite{Vak}.

Hodge integrals are intersection numbers of the form
    $$\langle\tau_{d_1}\cdots\tau_{d_n}\kappa_{a_1}\cdots\kappa_{a_m}\mid\lambda_{1}^{k_{1}}
    \cdots\lambda_{g}^{k_{g}}\rangle:=\int_{\overline{\mathcal{M}}_{g,n}}\psi_{1}^{d_{1}}\cdots\psi_{n}^{d_{n}}\kappa_{a_1}\cdots\kappa_{a_m}\lambda_{1}^{k_{1}}
    \cdots\lambda_{g}^{k_{g}},$$
which are rational numbers because the moduli space of curves are
orbifolds. They are nonzero only when
$\sum_{i=1}^nd_i+\sum_{i=1}^ma_i+\sum_{i=1}^gik_i=3g-3+n$.

Hodge integrals arise naturally in the localization computation of
Gromov-Witten invariants. They are extensively studied by
mathematicians and physicists. Hodge integrals involving only $\psi$
classes can be computed recursively by the the celebrated
Witten-Kontsevich theorem \cite{Wi, Ko}, which can be equivalently
formulated by the following DVV recursion relation \cite{DVV}
\begin{multline}\label{DVV}
\langle\tau_{k+1}\tau_{d_1}\cdots\tau_{d_n}\rangle_g=\frac{1}{(2k+3)!!}\left[\sum_{j=1}^n
\frac{(2k+2d_j+1)!!}{(2d_j-1)!!}\langle\tau_{d_1}\cdots
\tau_{d_{j}+k}\cdots\tau_{d_n}\rangle_g\right.\\
+\frac{1}{2}\sum_{r+s=k-1}
(2r+1)!!(2s+1)!!\langle\tau_r\tau_s\tau_{d_1}\cdots\tau_{d_n}\rangle_{g-1}\\
\left.+\frac{1}{2}\sum_{r+s=k-1} (2r+1)!!(2s+1)!!
\sum_{\underline{n}=I\coprod J}\langle\tau_r\prod_{i\in
I}\tau_{d_i}\rangle_{g'}\langle\tau_s\prod_{i\in
J}\tau_{d_i}\rangle_{g-g'}\right]
\end{multline}
where $\underline{n}=\{1,2,\ldots,n\}$.

Now there are several new proofs of Witten's conjecture \cite{CLL, Ka,
Ka-La, Ki-Li, Mi, Ok-Pa}.

Let {\it denom}$(r)$ denote the denominator of a rational number $r$
in reduced form (coprime numerator and denominator, positive
denominator). For $2g-2+n\geq1$, we define
$$D_{g,n}=lcm\left\{denom\left(\int_{\overline{\sM}_{g,n}}\psi_1^{d_1}\cdots\psi_n^{d_n}\right)\Big{|}\
\sum_{i=1}^{n}d_i=3g-3+n\right\},$$
and for $g\geq 2$,
$$\mathcal D_g=lcm\left\{denom\left(\int_{\overline{\sM}_{g,n}}\psi_1^{d_1}\cdots\psi_n^{d_n}\right)\Big{|}\
\sum_{i=1}^{n}d_i=3g-3+n,\ d_i\geq2, n\geq1\right\},$$
$$\widetilde{\mathcal D}_g=lcm\left\{denom\left(\int_{\overline{\sM}_{g}}\kappa_{a_1}\cdots\kappa_{a_m}\right)\Big{|}\
\sum_{i=1}^{m}a_m=3g-3\right\},$$
where {\it lcm} denotes the {\it
least common multiple}.

Note that $\mathcal D_g$ was previously defined by Itzykson and
Zuber \cite{It-Zu}.

We know that a neighborhood of $\Sigma\in\overline{\sM}_{g,n}$ is of
the form $U/\text{Aut}(\Sigma)$, where $U$ is an open subset of
$\mathbb C^{3g-3+n}$. This gives the orbifold structure for
$\overline{\sM}_{g,n}$. Since denominators of intersection numbers
on $\overline{\sM}_{g,n}$ all come from these orbifold quotient
singularities, the divisibility properties of $D_{g,n}$ and
$\mathcal D_g$ should reflect the overall behavior of singularities.

In Section 2, we study basic relations between $D_{g,n}$, $\mathcal
D_g$ and $\widetilde{\mathcal D}_g$. In Section 3, we discuss
briefly automorphism groups of Riemann surfaces and stable curves.
In Section 4, we study prime factors of $\mathcal D_g$ and prove a
strong form of a conjecture of Itzykson and Zuber \cite{It-Zu}
concerning denominators of intersection numbers. In Section 5, we
present a conjectural multinomial type property for intersection numbers
and verify it in low genera. \\

{\bf Acknowledgements.} We would like to thank Professor Sergei
Lando for valuable comments and interests in this work. The second
author is grateful to Professor Enrico Arbarello, Carel Faber, Sean
Keel, Rahul Pandharipande and Ravi Vakil for answering several
questions on moduli spaces of curves. We also thank the referees for
very helpful comments.

 \vskip
30pt
\section{Basic properties of $\mathcal D_g$}
If we take $k=-1$ and $k=0$ respectively in DVV formula~(\ref{DVV}), we get the string equation
$$\langle\tau_0\prod_{i=1}^{n}\tau_{k_i}\rangle_g=\sum_{j=1}^{n}\langle\tau_{k_j-1}\prod_{i\neq
j}\tau_{k_i}\rangle_g$$ and the dilaton equation
$$\langle\tau_1\prod_{i=1}^{n}\tau_{k_i}\rangle_g=(2g-2+n)\langle\prod_{i=1}^{n}\tau_{k_i}\rangle_g$$
Their proof may be found in the book \cite{La-Zv}.

\begin{lemma}
If $n\geq 1$, then
\begin{itemize}
\item[i)] $D_{0,n}=1$,

\item[ii)] $D_{1,n}=24$,

\item[iii)] $D_{g,1}=24^g\cdot g!$.
\end{itemize}

\end{lemma}
\begin{proof}
The lemma follows from the string equation, the dilaton equation and the following well known formulae
$$\langle\tau_{d_1}\cdots\tau_{d_n}\rangle_0=\binom{n-3}{d_1\cdots d_n}=\frac{(n-3)!}{d_1!\cdots d_n!},$$
$$\langle\tau_1\rangle_1=\frac{1}{24},\qquad \langle\tau_{3g-2}\rangle_g=\frac{1}{24^gg!}.$$
Their proofs can be found in \cite{La-Zv, Wi}.
\end{proof}

Note that $D_{0,n}=1$ is expected since $\overline{\sM}_{0,n}$ is a
smooth manifold.

\begin{theorem}
We have
$$D_{g,n}\mid D_{g,n+1}.$$
\end{theorem}

\begin{proof}
Let $q^s\mid D_{g,n}$, where $q$ is a prime number and $q^{s
+1}\nmid D_{g,n}$.

We sort $\{\langle\tau_{d_1}\dots\tau_{d_n}\rangle_g\mid
\sum_{i=1}^{n}d_i=3g-3+n,0\leq d_1\leq\dots\leq d_n\}$ in
lexicographical order, we say
$\langle\tau_{k_1}\dots\tau_{k_n}\rangle_g\prec\langle\tau_{m_1}\dots\tau_{m_n}\rangle_g$,
if there is some $i$, such that $k_j=m_j, j<i$ and $k_i<m_i$.

Let $\langle\tau_{k_1}\dots\tau_{k_n}\rangle_g$ be the minimal
element with respect to the lexicographical order such that its
denominator is divisible by $q^s$.

There exist integers $c,d,a_i,b_i$ where $i=1,\cdots,n-1$ such that
\begin{align*}
\langle\tau_0\tau_{k_1}\dots\tau_{k_n+1}\rangle_g &=
\langle\tau_{k_1}\dots\tau_{k_n}\rangle_g+
\sum_{i=1}^{n-1}\langle\tau_{k_1}\dots\tau_{k_i-1}\dots\tau_{k_{n-1}}\tau_{k_n+1}\rangle_g\\
&=\frac{c}{q^sd}+\sum_{i=1}^{n-1}\frac{b_i}{a_i},
\end{align*}
we require $q\nmid c, q\nmid d$ and $(a_i,b_i)=1$.

Since for $i=1,\dots,n-1$, we have\
$\langle\tau_{k_1}\dots\tau_{k_i-1}\dots\tau_{k_{n-1}}\tau_{k_n+1}\rangle_g\prec\langle\tau_{k_1}\dots\tau_{k_n}\rangle_g$,
so $a_i=q^{s_i}e_i$, where $s_i<l$ and $q\nmid e_i$. We now have

\begin{eqnarray*}
\langle\tau_0\tau_{k_1}\dots\tau_{k_n+1}\rangle_g=\frac{c\prod_{i=1}^{n-1}e_i+qd(\sum_{j=1}^{n-1}q^{s-s_j-1}\prod_{i\neq
j}e_i)}{q^sd\prod_{i=1}^{n-1}e_i}
\end{eqnarray*}
we see that $q$ can not divide the numerator, so we have proved
$q^s\mid D_{g,n+1}$. Since $q$ is arbitrary, we proved the
theorem.
\end{proof}

\begin{theorem}
We have $D_{g,n} \mid \widetilde{\mathcal D}_g$ for all $g\geq 2,
n\geq 1$. Moreover $\widetilde{\mathcal D}_g=D_{g,3g-3}$.
\end{theorem}

\begin{proof}

Let
$$
\pi_{n}:\overline{\sM}_{g,n}\longrightarrow
\overline{\sM}_{g,n-1},
$$
be the morphism that forgets the last marked point, then we have
 \cite{Ar-Co},
\begin{equation}
(\pi_{1}\dots\pi_n)_*(\psi_1^{a_1+1}\dots\psi_n^{a_n+1})=\sum_{\sigma\in
S_n}\kappa_\sigma,
\end{equation}
where $\kappa_\sigma$ is defined as follows. Write the permutation
$\sigma$ as a product of $\nu(\sigma)$ disjoint cycles, including
1-cycles: $\sigma=\beta_1\cdots\beta_{\nu(\sigma)}$, where we
think of the symmetric group $S_n$ as acting on the $n$-tuple
$(a_1,\dots ,a_n)$. Denote by $|\beta|$ the sum of the elements of
a cycle $\beta$. Then
$$
\kappa_\sigma=\kappa_{|\beta_1|}\kappa_{|\beta_2|}\dots
\kappa_{|\beta_{\nu(\sigma)}|}.
$$
From the formula (2), we get
$$\int_{\overline{\sM}_{g,n}}\psi_1^{a_1+1}\cdots\psi_n^{a_n+1}=\sum_{\sigma\in
S_n}\int_{\overline{\sM}_{g}}\kappa_\sigma,$$ so we proved $D_{g,n}
\mid \widetilde{\mathcal D}_g$.

On the other hand, any
$\int_{\overline{\sM}_{g}}\kappa_{a_1}\cdots\kappa_{a_m}$ can be
written as a sum of
$\int_{\overline{\sM}_{g,n}}\psi_1^{d_1}\cdots\psi_n^{d_n}$'s.
This can be seen by induction on the number of kappa classes, for
integrals with only one kappa class, we have
$\int_{\overline{\sM}_{g,n}}\kappa_{a_1}\psi_1^{d_1}\cdots\psi_n^{d_n}=\int_{\overline{\sM}_{g,n+1}}\psi_{n+1}^{a_1+1}\psi_1^{d_1}\cdots\psi_n^{d_n}$.
We also have
\begin{eqnarray*}
\int_{\overline{\sM}_{g,n}}\kappa_{a_1}\cdots\kappa_{a_m}\psi_1^{d_1}\cdots\psi_n^{d_n}&=&\int_{\overline{\sM}_{g,n+m}}\psi_{n+1}^{a_1+1}\cdots\psi_{n+m}^{a_m+1}\psi_1^{d_1}\cdots\psi_n^{d_n}\\
&&-\ \{\text{integrals with at most $m-1$ $\kappa$ classes}\}.
\end{eqnarray*}
thus finishing the induction argument. So we proved
$\widetilde{\mathcal D}_g=D_{g,3g-3}$.
\end{proof}

\begin{corollary}
For $g\geq 2$, we have $\mathcal D_g=\widetilde{\mathcal D}_g$.
\end{corollary}

We have computed $\mathcal D_g$ for $g\leq 20$ using the DVV formula
(\ref{DVV}) and observed the following conjectural exact values of $\mathcal
D_g$ (see also \cite{LX}).
\begin{conjecture}
Let $p$ be a prime number and $g\geq 2$. Let ${\rm ord}(p,n)$ denote
the maximum integer such that $p^{{\rm ord}(p,n)}\mid n$, then
\begin{itemize}
\item[i)] ${\rm ord}(2,\mathcal D_g)=3g+{\rm ord}(2,g!)$,

\item[ii)] ${\rm ord}(3,\mathcal D_g)=g+{\rm ord}(3,g!)$,

\item[iii)] ${\rm ord}(p,\mathcal D_g)=\lfloor\frac{2g}{p-1}\rfloor$ for
$p\geq5$, where $\lfloor x\rfloor$ denotes the maximum integer that
is not larger than $x$.
\end{itemize}
\end{conjecture}

On the other hand, we may get explicit expressions for multiples of
$\mathcal D_g$ by applying either Kazarian-Lando's formula
\cite{Ka-La} expressing intersection indices by Hurwitz numbers or
Proposition 4.4.

 \vskip
30pt
\section{Automorphism groups of stable curves}

First we recall some facts about automorphisms of compact Riemann surfaces following \cite{Fa-Kr}.

Let $X$ be a compact Riemann surface of genus $g$ and $Aut(X)$ the
group of conformal automorphisms of $X$. It's a classical theorem of
Hurwitz that if $g\geq 2$, then $|Aut(X)|\leq 84(g-1)$.

Let $G\subset Aut(X)$ be a group of automorphisms of $X$,
consider the natural map
$$\pi: X\rightarrow X/G$$
we know that $\pi$ has degree $|G|$ and $X/G$ is a compact Riemann
surface of genus $g_0$.

The mapping $\pi$ is branched only at the fixed points of $G$ and
the branching order
$$b(P)=\text{ord}G_P-1$$
where $G_P$ is the isotropy group at $P\in X$ which is known to be
cyclic.

Let $P_1,\ldots,P_r$ be a maximal set of inequivalent fixed points
of elements of $G\setminus\{1\}$. (that is, $P_i\neq h(P_j)$ for
all $h\in G$ and all $i\neq j$.)

Let $n_i=\text{ord}G_{P_i}$, then the total branch number of $\pi$
is given by

$$B=\sum_{i=1}^r\frac{|G|}{n_i}(n_i-1)=|G|\sum_{i=1}^{r}(1-\frac{1}{n_i})$$
the Riemann-Hurwitz formula now reads
\begin{equation*}
2g-2=|G|\left[2g_0-2+\sum_{i=1}^{r}(1-\frac{1}{n_i})\right]
\end{equation*}
so we have
\begin{equation}
|G|\ \Big|\ (2g-2)\cdot lcm(n_1,\ldots,n_r),
\end{equation}
this fact is crucial in the study of automorphism groups of
compact Riemann surfaces.

The following is a special case of a theorem due to W. Harvey (Theorem 6 in \cite{Ha}).
\begin{proposition}\cite{Ha}
The minimum genus $g$ of a compact Riemann surface which admits an
automorphism of order $p^r$ ($p$ is prime) is given by
$$g=\max\left\{2,\frac{p-1}{2}p^{r-1}\right\}.$$
\end{proposition}

In (3), we have $n_i=\text{ord}G_{P_i}$ and $G_{P_i}$ is cyclic, so Proposition 3.1 implies the following
\begin{corollary}
Let $X$ be a compact Riemann surface of genus $g\geq2$ and $G=|\text{Aut}(X)|$. Then
$${\rm ord}(p,|G|)\leq \lfloor\log_p\frac{2pg}{p-1}\rfloor+{\rm ord}(p,2(g-1)).$$
In particular, $p\nmid |G|$ if $p>2g+1$.
\end{corollary}

\begin{definition} A {\it node} on a curve is a point that is locally analytically isomorphic
to a neighborhood of the orgin of $xy=0$ in the complex plane
$\mathbb C^2$.

If $\Sigma$ is a nodal curve, define its {\it normalization}
$\tilde\Sigma$ to be the Riemann surface obtained by ``ungluing''
its nodes. Let $p: \tilde\Sigma\rightarrow\Sigma$ denote the
canonical normalization map. The preimages in $\tilde\Sigma$ of the
nodes of $\Sigma$ are called {\it node-branches}.

 A {\it stable curve} is a connected and compact nodal curve,
which means that its singular points are nodes and satisfy the
stability conditions: (i) each genus $0$ component has at least
three node-branches; (ii) each genus $1$ component has at least one
node-branch.

\end{definition}

Stability is equivalent to the finiteness of the automorphism group. Suppose $\Sigma$ is a stable curve of
arithmetic genus $g$ such that its normalization has $m$ components
$\Sigma_1,\dots,\Sigma_m$ of genus $g_1,\ldots,g_m$.

\begin{definition}
An automorphism $\varphi$ of the dual graph $\Gamma$ of $\Sigma$
will be called {\it geometric}, if it is induced by an automorphism
of the corresponding stable curve $\Sigma$. All geometric
automorphisms of $\Gamma$ form a group $GAut(\Gamma)$, which is a
subgroup of $Aut(\Gamma)$.
\end{definition}

The notion of geometric automorphism is introduced by Opstall and Veliche \cite{Op-Ve}
in their study of sharp bounds for the automorphism group of stable curves of a given genus.

\begin{theorem}
Let $\widetilde{Aut}(\Sigma_i)$ be the group of automorphisms of
$\Sigma_i$ fixing node-branches on $\Sigma_i$. Then we have
$$|Aut(\Sigma)|=|GAut(\Gamma)|\cdot\prod_{i=1}^{m}|\widetilde{Aut}(\Sigma_i)|$$
\end{theorem}
\begin{proof}
First note the following fact, if $f(x)$ and $g(y)$ are two
holomorphic functions defined near the origin of $\mathbb C^1$ and
satisfy $f(0)=g(0)$, then $F(x,y)=f(x)+g(y)-f(0)$ is a holomorphic
function near the origin of $\mathbb C^2$ satisfying $F(x,0)=f(x)$
and $F(0,y)=g(y)$. So to check whether a function on a nodal curve
is analytic, we need only check it is analytic restricting to each
connected component.

There is a natural map $p: Aut(\Sigma)\rightarrow GAut(\Gamma)$
mapping an automorphism of $\Sigma$ to the induced automorphism on
its dual graph $\Gamma$.

For each $b\in GAut(\Gamma)$, fix a $T_b\in Aut(\Sigma)$ such that
$p(T_b)=b$. If $f_i\in \widetilde{Aut}(\Sigma_i), i=1\cdots m$, we
denote by $(f_1,\cdots,f_m)\in Aut(\Sigma)$ the gluing morphism. We
define the following map
\begin{align*}
GAut(\Gamma)\times
\prod_{i=1}^{m}\widetilde{Aut}(\Sigma_i)&\longrightarrow Aut(\Sigma)\\
(b,f_1,\cdots,f_m)&\longmapsto T_b\circ (f_1,\cdots,f_m).
\end{align*}
It's not difficult to see that this map is in fact bijective. Its converse is
\begin{align*}
Aut(\Sigma)&\longrightarrow GAut(\Gamma)\times
\prod_{i=1}^{m}\widetilde{Aut}(\Sigma_i)\\
T&\longmapsto (p(T),(T_{p(T)}^{-1}\circ T) \big{|}_{\Sigma_1},\cdots,(T_{p(T)}^{-1}\circ T) \big{|}_{\Sigma_m}).
\end{align*}
So we
proved the theorem.
\end{proof}

\begin{proposition}
Let $\Sigma$ be a stable curve of arithmetic genus $g\geq 2$, if a
prime number $p$ divides $|\text{Aut}(\Sigma)|$, then $p\leq
2g+1$.
\end{proposition}
\begin{proof}
Let's assume that there are $\delta$ nodes on $\Sigma$ and
$\delta_i$ node-branches on each $\Sigma_i$. Then we have the
following relations,
\begin{align}
&g=\sum_{i=1}^{m}(g_i-1)+\delta+1,\\
&2g_i+\delta_i-2\geq 1,\\
&2\delta=\sum_{i=1}^{m}\delta_i.
\end{align}
Sum up (5) for $i=1$ to $n$ and substitute (4) and (6) into (5),
we get
$$m\leq 2g-2.$$

Let $e_{ij}$ denote the number of edges between $\Sigma_i$ and
$\Sigma_j$ in the dual graph of $\Sigma$, then it's obvious that
$e_{ij}\leq g+1$.

Since $|\text{Aut}(\Gamma)|$ divides $m!\prod_{(i,j)}(e_{ij}!)$
which is not divisible by prime numbers greater than $2g+1$, and
$g_i\leq g$, so the proposition follows from Theorem 3.5 and
Corollary 3.2.
\end{proof}

We remark that for non-stable nodal curves, Proposition 3.6 may not
hold.

\vskip 30pt
\section{Prime factors of $\mathcal D_g$}

\begin{definition} In \cite{Fa2},
the following generating function
$$F(x_1,\cdots,x_n)=\sum_{g=0}^{\infty}\sum_{\sum d_i=3g-3+n}\langle\tau_{d_1}\cdots\tau_{d_n}\rangle_g\prod_{i=1}^n x_i^{d_i}$$
is called the $n$-point function.
\end{definition}

In particular, 2-point function has a simple explicit form due to
Dijkgraaf (see \cite{Fa2})
$$F(x_1,x_2)=\frac{1}{x_1+x_2}\exp\left(\frac{x_1^3}{24}+\frac{x_2^3}{24}\right)\sum_{k=0}^{\infty}\frac{k!}{(2k+1)!}\left(\frac{1}{2}x_1x_2(x_1+x_2)\right)^k.$$

\begin{lemma}
Let $p$ be a prime number and $g\geq 2$, then
\begin{itemize}
\item[i)] If $p>2g+1$, then $p\nmid D_{g,2}$,

\item[ii)] If $g+1\leq p\leq 2g+1$, then $$p\mid
denom\langle\tau_{\frac{p-1}{2}}\tau_{3g-1-\frac{p-1}{2}}\rangle_g,$$

\item[iii)] If $2g+1$ is prime, then $(2g+1)\mid
denom\langle\tau_{d}\tau_{3g-1-d}\rangle_g$ if and only if $g\leq
d\leq 2g-1$.

\item[iv)] If $2g+1$ is prime, then ${\rm ord}(2g+1,D_{g,2})=1$.
\end{itemize}
\end{lemma}
\begin{proof}
From the 2-point function, we get
\begin{align*}
\langle\tau_{d}\tau_{3g-1-d}\rangle_g=&\sum_{i=0}^g\sum_k\binom{g-k}{i}\binom{k-1}{d-3i-k}\frac{k!}{(g-k)!24^{g-k}(2k+1)!2^k}\\
&+\frac{(-1)^{d\bmod3}}{g!24^g}\binom{g-1}{\lfloor\frac{d}{3}\rfloor},
\end{align*}
where the summation range of $k$ is
$\max(\frac{d_1-3i+1}{2},1)\leq k\leq \min(g-i,d_1-3i)$. Then the
lemma follows easily.
\end{proof}

\begin{theorem}
Let $p$ be a prime number, $g\geq 2$ and let ${\rm ord}(p,q)$ denote
the maximum integer such that $p^{{\rm ord}(p,q)}\mid q$, then
\begin{itemize}
\item[i)] If $p>2g+1$, then $p\nmid \mathcal D_g$,

\item[ii)] For any prime $p\leq 2g+1$, we have $p\mid \mathcal D_g$,

\item[iii)] If $2g+1$ is prime, then ${\rm ord}(2g+1,\mathcal D_g)=1$,

\item[iv)] ${\rm ord}(2,\mathcal D_g)=3g+{\rm ord}(2,g!)$.
\end{itemize}

\end{theorem}
\begin{proof}
For part (i), we use induction on the pair of genus and the number
of marked points $(g,n)$ to prove that denominators of all $\psi$
class intersection numbers
$\langle\tau_{d_1}\cdots\tau_{d_n}\rangle_g$ are not divisible by
prime numbers greater than $2g+1$. If $p>2g+1$, then $p\nmid
D_{g,2}$ has been proved in Lemma 4.2(i). Also $\mathcal
D_2=2^7\cdot3^2\cdot5$ is not divisible by $p>5$. So we may assume
$g\geq3, n\geq3$. We rewrite the DVV formula here,
\begin{eqnarray*}
\langle\tau_{d_1}\cdots\tau_{d_n}\rangle_g&=&\frac{1}{(2d_1+1)!!}\left[\sum_{j=2}^n
\frac{(2d_1+2d_j-1)!!}{(2d_j-1)!!}\langle\tau_{d_2}\cdots
\tau_{d_{j}+d_1-1}\cdots\tau_{d_n}\rangle_g\right.\\
&&+\frac{1}{2}\sum_{r+s=d_1-1}
(2r+1)!!(2s+1)!!\langle\tau_r\tau_s\tau_{d_2}\cdots\tau_{d_n}\rangle_{g-1}\\
&&\left.+\frac{1}{2}\sum_{r+s=d_1-1} (2r+1)!!(2s+1)!!
\sum_{\{2,\cdots,n\}=I\coprod J}\langle\tau_r\prod_{i\in
I}\tau_{d_i}\rangle_{g'}\langle\tau_s\prod_{i\in
J}\tau_{d_i}\rangle_{g-g'}\right]
\end{eqnarray*}

For $n\geq 3$ marked points, we may take $d_1\leq g$, then by
induction on $(g,n)$ it's easy to see that the denominator of the
right hand side is not divisible by prime numbers greater than
$2g+1$.

For part (ii), it follows from Lemma 2.1(iii), Theorem 2.3 and Lemma
4.2(ii).

For part (iii), we again use induction on $(g,n)$ as in the proof of
part (i), we may assume $n\geq3$. In view of Lemma 4.2(iii)-(iv), we
need only prove ${\rm ord}(2g+1,D_{g,n})\leq1$. If $n>3$, then we
may take $d_1<g$ in $\langle\tau_{d_1}\cdots\tau_{d_n}\rangle_g$,
whose denominator is not divisible by $(2g+1)^2$. This is easily
seen by induction on the right hand side of the DVV formula. So we
are only left to prove that the denominator of
$\langle\tau_g\tau_g\tau_g\rangle_g$ is not divisible by $(2g+1)^2$.
\begin{eqnarray*}
\langle\tau_g\tau_g\tau_g\rangle_g=\frac{1}{(2g+1)!!}\left[\frac{2(4g-1)!!}{(2g-1)!!}
\langle\tau_g\tau_{2g-1}\rangle_g+\{\text{lower genus
terms}\}\right]
\end{eqnarray*}
Since the factor $2g+1$ in the denominator of
$\langle\tau_g\tau_{2g-1}\rangle$ will be cancelled by $(4g-1)!!$,
by induction we proved (iii).

For part (iv), since $\langle\tau_{3g-2}\rangle_g=\frac{1}{24^gg!}$,
we have ${\rm ord}(2,\mathcal D_g)\geq 3g+{\rm ord}(2,g!)$, the
reverse inequality can be seen from the DVV formula by
induction on $(g,n)$ and note the following,
\begin{eqnarray*}
&&\frac{1}{2}\sum_{r+s=k-1} (2r+1)!!(2s+1)!!
\sum_{\underline{n}=I\coprod J}\langle\tau_r\prod_{i\in
I}\tau_{d_i}\rangle_{g'} \langle\tau_s\prod_{i\in
J}\tau_{d_i}\rangle_{g-g'}\\&=& \sum_{r+s=k-1} (2r+1)!!(2s+1)!!
\sum_{\{2,\dots,n\}=I\coprod J}\langle\tau_r\tau_{d_1}\prod_{i\in
I}\tau_{d_i}\rangle_{g'} \langle\tau_s\prod_{i\in
J}\tau_{d_i}\rangle_{g-g'}.
\end{eqnarray*}
\end{proof}

\begin{lemma}
If $2\leq p\leq g+1$ is a prime number, then ${\rm ord}(p,
D_{g,3})\geq 2$.
\end{lemma}
\begin{proof} From Lemma 2.1(3), we have $24^g\mid D_{g,3}$, so the lemma is obvious for $p=2$ or
$3$. We assume $p\geq 5$ below.

The following formula of the special three-point function is due to
Faber \cite{Fa2}.
\begin{align*}
F_g(x,y,-y)&=\sum_{b\geq0}\sum_{j=0}^{2b}(-1)^j\langle\tau_{3g-2b}\tau_j\tau_{2b-j}\rangle_g x^{3g-2b}y^{2b}\\
&=\sum_{\substack{a+b+c=g\\b\geq a}}\frac{(a+b)!}{4^{a+b}24^c
(2a+2b+1)!!(b-a)!(2a+1)!c!}x^{3a+3c+b}y^{2b}.
\end{align*}

If $p>\frac{2g+1}{3}$, then consider the coefficient of
$x^{3g-p+1}y^{p-1}$ in $F_g(x,y,-y)$
$$[F_g(x,y,-y)]_{x^{3g-p+1}y^{p-1}}=\sum_{\substack{a+b+c=g\\a\leq \frac{p-1}{2}}}\frac{(a+b)!}{4^{a+b}24^c
(2a+2b+1)!!(b-a)!(2a+1)!c!}$$ where $b=\frac{p-1}{2}$. We must have
$c<p$, so it's not difficult to see that only the term with
$a=b=\frac{p-1}{2}$ can contain factor $p^2$ in the denominator.

If $p\leq \frac{2g+1}{3}$, we have
$$[F_g(x,y,-y)]_{x^gy^{2g}}=\frac{1}{4^g(2g+1)!!},$$
and ${\rm ord}(p,(2g+1)!!)\geq 2$.

So we proved the lemma.
\end{proof}

\begin{theorem}
Let $X$ be a compact Riemann surface of genus $g'\geq 2$ and $g\geq
g'$, then $|Aut(X)|$ divides $D_{g,3}$.
\end{theorem}
\begin{proof} We first prove the case $g'=g$.

Let $p$ denote a prime number. By Corollary 3.2, it is sufficient to
prove
\begin{equation}
\lfloor\log_p\frac{2pg}{p-1}\rfloor+{\rm ord}(p,2(g-1))\leq {\rm
ord}(p,D_{g,3})
\end{equation}
for all prime $p\leq 2g+1$.

If $\max(g,5)\leq p\leq 2g+1$, then we have
$\lfloor\log_p\frac{2pg}{p-1}\rfloor\leq 1$ and ${\rm
ord}(p,2(g-1))=0$, so from Theorem 4.3(2), the above inequality (8)
holds in this case.

Now we assume $5\leq p\leq g-1$, the cases $p=2$ and $p=3$ will be
treated at last. We still need to divide into three finer cases.

{\it Case i)}\  If $p=g-1\geq 5$ is prime, then we have
$(g-1)(g-2)>2g$. By Lemma 4.7, we have
$$\lfloor\log_{g-1}\frac{2g(g-1)}{g-2}\rfloor+1\leq 2\leq {\rm ord}(g-1,D_{g,3}).$$

{\it Case ii)}\  Otherwise if $p\nmid(g-1)$, since ${\rm
ord}(p,2(g-1))=0$, $g!\mid D_{g,3}$ and ${\rm
ord}(p,g!)\geq\lfloor\frac{g}{p}\rfloor$, so in order to check (8),
it's sufficient to prove
$$\lfloor\log_p\frac{2pg}{p-1}\rfloor\leq\lfloor\frac{g}{p}\rfloor.$$

Let $g=kp+r$, where $\ -p\leq r<0$. Then
$\lfloor\frac{g}{p}\rfloor=k-1$. Since for fixed $k$, the left hand
side takes its maximum value when $g=kp-1$, we need only prove the
above identity for $g=kp-1$, which is equivalent to for all $k\geq
2$, $p\geq5$,
\begin{equation*}
p^k>\frac{2p(kp-1)}{p-1}, \quad i.e.\ \ p^{k}-p^{k-1}-2kp+2>0,
\end{equation*}
which is not difficult to check.

{\it Case iii)}\  If $p\mid(g-1)$ and $5\leq p<g-1$. Let ${\rm
ord}(p,2(g-1))=r$. Then $p^r\mid (g-1)$, we have
\begin{align*}
{\rm ord}(p,D_{g,3})\geq{\rm
ord}(p,g!)&=\lfloor\frac{g}{p}\rfloor+\lfloor\frac{g}{p^2}\rfloor
+\lfloor\frac{g}{p^3}\rfloor+\cdots\\
&\geq\lfloor\frac{g}{p}\rfloor+r-1.
\end{align*}

So it's sufficient to prove
$$\lfloor\log_p\frac{2pg}{p-1}\rfloor+1\leq\lfloor\frac{g}{p}\rfloor.$$

Let $g=kp+1,\ k\geq2$, we need to prove
\begin{equation*}
p^k>\frac{2p(kp+1)}{p-1}, \quad i.e.\ \ p^{k}-p^{k-1}-2kp-2>0.
\end{equation*}
The above inequality holds except in the case $p=5$, $k=2$ and
$g=11$, which should be treated separately. We have
$${\rm ord}(5,|G|)\leq \lfloor\log_5\frac{110}{4}\rfloor+1=3$$
and ${\rm ord}(5, D_{11,3})=3$, in fact
$$D_{11,3}=2^{41}\cdot 3^{15}\cdot 5^3\cdot 7^2\cdot 11^2\cdot 13\cdot 17\cdot 19\cdot 23.$$
We finished checking in this case.

Now we consider the remaining two cases, $p=2$ and $p=3$. Note that
$24^gg!\mid D_{g,3}$.

If $p=2$, it's sufficient to prove $\log_2{4g}\leq 3g-1$.

If $p=3$, it's sufficient to prove $\log_3{3g}\leq g$.

Both cases are easy to check. So we conclude the proof of the
theorem when $g'=g$.

The proof of the cases $g'<g$ can be proved by exactly the same
argument and using Lemma 4.7.
\end{proof}

We remark that there exists a compact Riemann surface $X$ of genus 6 with
$|Aut(X)|=150$ (see Table 13 in \cite{Br}). While the power of $5$ in $D_{6,2}=2^{22}\cdot 3^8\cdot 5\cdot 7\cdot 11\cdot 13$
is only $1$, so $|Aut(X)|\nmid D_{6,2}$.
In this sense, we may say that Theorem 4.8 is optimal.

The following immediate corollary of Theorem 4.8 is a conjecture of
Itzykson and Zuber, stated at the end of Section 5 of \cite{It-Zu}.
\begin{corollary}
For $1<g'\le g$, the order of automorphism group of any compact
Riemann surface of genus $g'$ divides $\mathcal D_g$.
\end{corollary}

We remark that the statement of Corollary 4.9 doesn't hold for
stable curves, namely there exists some stable curve of genus $g$,
the order of whose automorphism group does not divide $\mathcal
D_g$. A counterexample can be constructed as follows. Let $n=\lfloor
\frac{2g}{p-1}\rfloor$ Riemann surfaces of genus $\frac{p-1}{2}$
attached to a sphere at $e^{\frac{2\pi i}{n}}$ for $0\leq i\leq
n-1$. When $n\geq p$, the order of automorphism group of such a
stable curve will have the power of $p$ larger than $\lfloor
\frac{2g}{p-1}\rfloor$ (see conjecture 2.5).

\vskip 30pt
\section{A conjectural numerical property of intersection numbers}

During our work on intersection numbers, we noticed a multinomial type property
for intersection numbers. Although still conjectural, we feel they are interesting constraints of intersection numbers on
moduli spaces, so we briefly present them here.

From
$$\langle\tau_{d_1}\cdots\tau_{d_n}\rangle_0=\binom{n-3}{d_1\cdots d_n}=\frac{(n-3)!}{d_1!\cdots d_n!},$$
we see that if $d_1<d_2$, we have
$$\langle\tau_{d_1}\tau_{d_2}\cdots\tau_{d_n}\rangle_0\leq
\langle\tau_{d_1+1}\tau_{d_2-1}\cdots\tau_{d_n}\rangle_0.$$

Now we prove that the same inequality holds in genus $1$.
\begin{proposition}
For $\sum_{i=1}^{n}d_i=n$ and $d_1<d_2$, we have
$$\langle\tau_{d_1}\tau_{d_2}\cdots\tau_{d_n}\rangle_1\leq
\langle\tau_{d_1+1}\tau_{d_2-1}\cdots\tau_{d_n}\rangle_1.$$
\end{proposition}
\begin{proof}
We prove the inequality by induction on $n$. If $n=2$, we have
$$\langle\tau_0\tau_2\rangle_1=\langle\tau_1\tau_1\rangle_1=\frac{1}{24}.$$

Now assume that the proposition has been proved for $n-1$. We may
also assume $d_2-d_1\geq 2$, otherwise it is trivial. So by the
symmetry property of intersection numbers, we may assume without
loss of generality that $d_n=0$ or $d_n=1$.

If $d_n=1$ then by dilaton equation
\begin{eqnarray*}
\langle\tau_{d_1}\tau_{d_2}\cdots\tau_{d_n}\rangle_1&=&(n-1)\langle\tau_{d_1}\tau_{d_2}\cdots\tau_{d_{n-1}}\rangle_1\\
\langle\tau_{d_1+1}\tau_{d_2-1}\cdots\tau_{d_n}\rangle_1&=&(n-1)\langle\tau_{d_1+1}\tau_{d_2-1}\cdots\tau_{d_{n-1}}\rangle_1.
\end{eqnarray*}

So $\langle\tau_{d_1}\tau_{d_2}\cdots\tau_{d_n}\rangle_1\leq
\langle\tau_{d_1+1}\tau_{d_2-1}\cdots\tau_{d_n}\rangle_1$ holds in
this case by induction.

If $d_n=0$ then by string equation
\begin{eqnarray*}
\langle\tau_{d_1}\tau_{d_2}\cdots\tau_{d_n}\rangle_1&=&\langle\tau_{d_1-1}\tau_{d_2}\cdots\tau_{d_{n-1}}\rangle_1+
\langle\tau_{d_1}\tau_{d_2-1}\cdots\tau_{d_{n-1}}\rangle_1\\
&&+\sum_{i=3}^{n-1}\langle\tau_{d_1}\tau_{d_2}\cdots\tau_{d_i-1}\cdots\tau_{d_{n-1}}\rangle_1\\
\langle\tau_{d_1+1}\tau_{d_2-1}\cdots\tau_{d_n}\rangle_1&=&\langle\tau_{d_1}\tau_{d_2-1}\cdots\tau_{d_{n-1}}\rangle_1+
\langle\tau_{d_1+1}\tau_{d_2-2}\cdots\tau_{d_{n-1}}\rangle_1\\
&&+\sum_{i=3}^{n-1}\langle\tau_{d_1+1}\tau_{d_2-1}\cdots\tau_{d_i-1}\cdots\tau_{d_{n-1}}\rangle_1.
\end{eqnarray*}

So $\langle\tau_{d_1}\tau_{d_2}\cdots\tau_{d_n}\rangle_1\leq
\langle\tau_{d_1+1}\tau_{d_2-1}\cdots\tau_{d_n}\rangle_1$ holds
again by induction.
\end{proof}

Now we formulate the following conjecture
\begin{conjecture}
For $\sum_{i=1}^{n}d_i=3g-3+n$ and $d_1<d_2$, we have
$$\langle\tau_{d_1}\tau_{d_2}\cdots\tau_{d_n}\rangle_g\leq
\langle\tau_{d_1+1}\tau_{d_2-1}\cdots\tau_{d_n}\rangle_g.$$
\end{conjecture}
Namely the more evenly $3g-3+n$ be distributed among indices, the
larger the intersection numbers.

By the same argument of Proposition 5.1, we can see that for each
$g$, it's enough to check only those intersection numbers with
$n\leq 3g-1$ and $d_3\geq 2,\ldots,d_n\geq 2$.

We checked Conjecture 5.2 for $g\leq 16$ with the help of Faber's
Maple program. Moreover, for $n=2$, we checked all $g\leq 300$ using
Dijkgraaf's $2$-point function; for $n=3$, we checked all $g\leq 50$
using Zagier's $3$-point function.

$$ \ \ \ \ $$
\bibliographystyle{amsplain}

\end{document}